\documentclass[11pt]{article}
\usepackage{graphicx}
\usepackage{amsfonts}
\usepackage{theorem}
\newtheorem{Lem}{Lemma}[section]

\newtheorem{The}[Lem]{Theorem}
\newtheorem{Prop}[Lem]{Proposition}
\newtheorem{Cor}[Lem]{Corollary}

\newtheorem{Prob}[Lem]{Problem}

\newtheorem{Con}[Lem]{Conjecture}
\newcommand{\qed}{\hbox{\rule{6pt}{6pt}}}

\begin{document}
\title{A matrix trace inequality and its application}
\author{Shigeru Furuichi$^1$\footnote{E-mail:furuichi@chs.nihon-u.ac.jp} and Minghua Lin$^2$\footnote{E-mail:lin243@uregina.ca}\\
$^1${\small Department of Computer Science and System Analysis,}\\
{\small College of Humanities and Sciences, Nihon University,}\\
{\small 3-25-40, Sakurajyousui, Setagaya-ku, Tokyo, 156-8550, Japan}\\
$^2${\small Department of Mathematics and Statistics,}\\
{\small University of Regina, Regina, Saskatchewan, S4S 0A2, Canada}}
\date{}
\maketitle
{\bf Abstract.} In this short paper, we give a complete and affirmative answer to a conjecture on 
matrix trace inequalities for the sum of positive semidefinite matrices.
We also apply the obtained inequality to derive a kind of generalized Golden-Thompson inequality for
positive semidefinite matrices.
\vspace{3mm}


\vspace{3mm}
\vspace{3mm}

\section{Introduction}
We give some notations.
The set of all $n \times n$ matrices on the complex field $\mathbb{C}$ is represented by $M(n,\mathbb{C})$.
The set of all $n \times n$ Hermitian matrices is also represented by $M_h(n,\mathbb{C})$.
Moreover the set of all $n \times n$ nonnegative (positive semidefinite) matrices is also represented by $M_+(n,\mathbb{C})$.
Here $X\in M_+(n,\mathbb{C})$ means we have $\langle \phi \vert X \vert \phi \rangle \geq 0$ for any vector
$\vert \phi \rangle \in \mathbb{C}^n$.

The purpose of this short paper is to give the answer to the following conjecture which was given in the paper \cite{Furu1}.
\begin{Con}  {\bf (\cite{Furu1})}\label{con}
For $X,Y\in M_+(n,\mathbb{C})$ and $p\in \mathbb{R}$, the following inequalities hold or not?
\begin{itemize}
\item[(i)] $Tr[(I + X +Y +Y^{1/2}XY^{1/2})^p] \leq Tr[(I+X+Y+XY)^p]$ for $p \geq 1$.
\item[(ii)] $Tr[(I + X +Y +Y^{1/2}XY^{1/2})^p] \geq Tr[(I+X+Y+XY)^p]$ for $0 \leq p \leq 1$.
\end{itemize}
\end{Con}

We firstly note that the matrix $I+X+Y+XY=(I+X)(I+Y)$ is generally not positive semidefinite.
However, the eigenvalues of the matrix $(I+X)(I+Y)$ are same to those of the positive semidefinite matrix
$(I+X)^{1/2}(I+Y)(I+X)^{1/2}$. Therefore the expression $Tr[(I+X+Y+XY)^p]$ always makes sense.

We easily find that the equality for (i) and (ii) in Conjecture \ref{con} holds in the case of  $p=1$. In addition,
the case of $p=2$ was proven by elementary calculations in \cite{Furu1}.

Putting $T=(I+X)^{1/2}$ and $S=Y^{1/2}$, Conjecture \ref{con} can be reformulated by the following problem,
because we have
$Tr[(I+X+Y+XY)^p]=Tr[(T^2+T^2S^2)^p]=Tr[(T^2(I+S^2))^p]=Tr[(T(I+S^2)T)^p]=Tr[(T^2+TS^2T)^p].$

\begin{Prob} \label{prob02}
For $T,S\in M_+(n,\mathbb{C})$ and $p\in \mathbb{R}$, the following inequalities hold or not?
\begin{itemize}
\item[(i)] $ Tr[(T^2+ST^2S)^p] \leq Tr[(T^2+TS^2T)^p]$ for  $p \geq 1$.
\item[(ii)] $ Tr[(T^2+ST^2S)^p] \geq Tr[(T^2+TS^2T)^p]$ for $0 \leq p \leq 1$.
\end{itemize}
\end{Prob}

\section{Main results}
To solve Problem \ref{prob02}, we use the concept of the majorization. See \cite{MO} for the details on the majorization.
Here for $X\in M_h(n,\mathbb{C})$,
$\lambda^{\downarrow }(X)=\left(\lambda_1^{\downarrow}  (X),\cdots,\lambda_n^{\downarrow}  (X) \right)$
 represents the eigenvalues of the Hermitian matrix $X$ in decreasing order,
$\lambda_1^{\downarrow}  (X) \geq \cdots \geq \lambda_n^{\downarrow}  (X)$.
In addition $x \prec y$ means that $x=(x_1,\cdots,x_n)$ is majorized by $y=(y_1,\cdots,y_n)$, if we have
$$\sum_{j=1}^kx_j \leq \sum_{j=1}^k y_j\quad (k=1,\cdots,n-1)$$
and
$$\sum_{j=1}^nx_j =\sum_{j=1}^n y_j.$$

We need the following lemma which can be obtained as a consequence of Ky Fan's maximum principle.
\begin{Lem} {\bf (p.35 in \cite{Bha})}   \label{Lind-lem}
For $A,B \in M_h(n,\mathbb{C})$ and any $k=1,2,\cdots,n$, we have
\begin{equation}
\sum_{j=1}^k \lambda_{j}^{\downarrow}(A+B) \leq \sum_{j=1}^k \lambda_{j}^{\downarrow}(A)+\sum_{j=1}^k \lambda_{j}^{\downarrow}(B).
\end{equation}
\end{Lem}

Then we have the following theorem.
\begin{The} \label{the-01}
For $S,T\in M_+(n,\mathbb{C})$, we have
\begin{equation}   \label{maj_ineq00}
\lambda^{\downarrow }(T^2+ST^2S) \prec \lambda^{\downarrow }(T^2+TS^2T)
\end{equation}
\end{The}

{\it Proof}:
For $S,T\in M_+(n,\mathbb{C})$, we need only to show the following
\begin{equation}   \label{maj_ineq01}
\sum_{j=1}^{k} \lambda_{j}^{\downarrow}  (T^2+ST^2S) \leq \sum_{j=1}^{k} \lambda_{j}^{\downarrow}   (T^2+TS^2T)
\end{equation}
for $k=1,2,\cdots,n-1$,
since we have
$$
\sum_{j=1}^n \lambda_j ^{\downarrow}  (T^2+ST^2S) =\sum_{j=1}^n \lambda_j ^{\downarrow}  (T^2+TS^2T),
$$
which is equivalent to $Tr[T^2+ST^2S]=Tr[T^2+TS^2T]$.

By  Lemma \ref{Lind-lem}, we have
\begin{equation} \label{Lind-ineq01}
2 \sum_{j=1}^k \lambda_{j}^{\downarrow}(X) \leq
\sum_{j=1}^k \lambda_{j}^{\downarrow}\left( X+Y \right)+\sum_{j=1}^k \lambda_{j}^{\downarrow}\left( X-Y \right).
\end{equation}
for $X,Y\in M_h(n,\mathbb{C})$ and any $k=1,2,\cdots,n$.

For $X\in M(n,\mathbb{C})$, the matrices $XX^*$ and $X^*X$ are unitarily similar so that we have $ \lambda_j^{\downarrow}(XX^*) =\lambda_j^{\downarrow}(X^*X) $.
Then we have the following inequality:
\begin{eqnarray*}
2 \sum_{j=1}^{k} \lambda_{j}^{\downarrow}  \left( T^2+TS^2T \right) &=&
\sum_{j=1}^{k} \lambda_{j}^{\downarrow}\left(  T^2+TS^2T \right)  +\sum_{j=1}^{k} \lambda_{j}^{\downarrow}\left( T^2+TS^2T \right)  \\
&=& \sum_{j=1}^{k} \lambda_{j}^{\downarrow}\left( (T+iTS)(T-iST)\right)  +\sum_{j=1}^{k} \lambda_{j}^{\downarrow}\left( (T-iTS)(T+iST)\right)  \\
&=& \sum_{j=1}^{k} \lambda_{j}^{\downarrow}\left(  (T-iST)(T+iTS)\right)  +\sum_{j=1}^{k} \lambda_{j}^{\downarrow}\left(  (T+iST)(T-iTS)\right)  \\
&=& \sum_{j=1}^{k} \lambda_{j}^{\downarrow}\left( T^2+ST^2S+i\left( T^2S-ST^2 \right)\right) +
\sum_{j=1}^{k} \lambda_{j}^{\downarrow}\left(  T^2+ST^2S-i\left( T^2S-ST^2 \right) \right)  \\
&\geq& 2 \sum_{j=1}^{k} \lambda_{j}^{\downarrow} \left( T^2+ST^2S\right),
\end{eqnarray*}
for any $k=1,2,\cdots,n-1$,
by using the inequality (\ref{Lind-ineq01}) for $X=T^2+ST^2S$ and $Y=i(T^2S-ST^2)$.
Thus we have the inequality (\ref{maj_ineq01}) so that the proof is completed.

\hfill \qed

From Theorem \ref{the-01}, we have the following corollary.

\begin{Cor} \label{the01}
For $T,S\in M_+(n,\mathbb{C})$ and $p\in \mathbb{R}$, the following inequalities hold.
\begin{itemize}
\item[(i)] $ Tr[(T^2+ST^2S)^p] \leq Tr[(T^2+TS^2T)^p]$ for  $p \geq 1$.
\item[(ii)] $ Tr[(T^2+ST^2S)^p] \geq Tr[(T^2+TS^2T)^p]$ for $0 \leq p \leq 1$.
\end{itemize}
\end{Cor}

{\it Proof :}
Since $f(x)=x^p,(p\geq 1)$ is convex function and $f(x)=x^p,(0\leq p \leq 1)$ is concave function,
we have the present corollary thanks to Theorem \ref{the-01} and a general property of majorization (See p.40 in \cite{Bha}).
\hfill \qed

As mentioned in Introduction, Corollary \ref{the01} implies the following corollary by putting
$T=(I+X)^{1/2}$ and $S=Y^{1/2}$.

\begin{Cor}  \label{con-cor}
For $X,Y\in M_+(n,\mathbb{C})$ and $p\in \mathbb{R}$, the following inequalities hold.
\begin{itemize}
\item[(i)] $Tr[(I + X +Y +Y^{1/2}XY^{1/2})^p] \leq Tr[(I+X+Y+XY)^p]$ for $p \geq 1$.
\item[(ii)] $Tr[(I + X +Y +Y^{1/2}XY^{1/2})^p] \geq Tr[(I+X+Y+XY)^p]$ for $0 \leq p \leq 1$.
\end{itemize}
\end{Cor}

Thus Conjecture \ref{con} was completely solved with an affirmative answer.


\section{An application}
In this section, we give a kind of one-parameter extension of the famous Golden-Thompson inequality \cite{Gol,Tho} for
positive semidefinite matrices,  applying the obtained result in the previous section.
For this purpose, we denote the generalized exponential function by $\exp_{\nu}(X) \equiv \left( I +\nu X \right) ^{\frac{1}{\nu}}$
for $\nu \in (0,1]$ and $X\in M(n,\mathbb{C})$ such that $Tr[(I +\nu X)^{\frac{1}{\nu}}] \in \mathbb{R}$.
In addition, we use the following inequalities proved in \cite{Furu2}.

\begin{Lem}  {\bf (\cite{Furu2})} \label{GT_ineq_gen_before}
For $X,Y\in M_+(n,\mathbb{C})$, and $\nu \in (0,1]$, we have
\begin{itemize}
\item[(i)]
\begin{equation}   \label{GT_ineq_gen_EQ00}
Tr[  \exp_{\nu}(X+Y)  ] \leq Tr[  \exp_{\nu}(X+Y+\nu Y^{1/2}XY^{1/2})  ].
\end{equation}
\item[(ii)]
\begin{equation}   \label{GT_ineq_gen_EQ}
Tr[\exp_{\nu}(X+Y+\nu XY)] \leq Tr[\exp_{\nu}(X) \exp_{\nu}(Y)].
\end{equation}
\end{itemize}
\end{Lem}

As mentioned in the below of Conjecture \ref{con}, the expression of the left hand side in (\ref{GT_ineq_gen_EQ}) makes also sense, 
since we have $Tr[\exp_{\nu}(X+Y+\nu XY)] = Tr[\left\{ (I+\nu X)^{1/2} (I+\nu Y) (I+\nu X)^{1/2}  \right\}^{\frac{1}{\nu}}] \geq 0$. 

From (i) of Corollary \ref{con-cor} and Lemma \ref{GT_ineq_gen_before}, we have the following proposition.

\begin{Prop}   \label{gen-GT}
For $X,Y\in M_+(n,\mathbb{C})$ and $\nu \in (0,1]$, we have
\begin{equation}\label{gen-GT-ineq}
Tr[  \exp_{\nu}(X+Y)  ] \leq Tr[\exp_{\nu}(X) \exp_{\nu}(Y)].
\end{equation}
\end{Prop}

{\it Proof}:
The right hand side of (\ref{GT_ineq_gen_EQ00}) is bounded from the above by applying (i) of Corollary \ref{con-cor} and
putting $X_1=\nu X$, $Y_1=\nu Y$ and $p=\frac{1}{\nu}$:
\begin{eqnarray*}
Tr\left[\exp_{\nu}(X+Y+\nu Y^{1/2}XY^{1/2})\right] &=& Tr\left[ \left\{ I+\nu (X+Y+\nu Y^{1/2}XY^{1/2}) \right\}^{\frac{1}{\nu}}\right] \\
&=& Tr \left[ (I+X_1+Y_1+Y_1^{1/2}X_1Y_1^{1/2})^p \right]\\
&\leq &  Tr \left[ (I+X_1+Y_1+X_1Y_1)^p \right]\\
&=& Tr \left[ \left\{ I+\nu (X+Y+\nu XY) \right\}^{\frac{1}{\nu}} \right]\\
&=& Tr\left[ \exp_{\nu}(X+Y+\nu XY) \right],
\end{eqnarray*}
which is the left hand side of (\ref{GT_ineq_gen_EQ}). Thus we have the present proposition
thanks to Lemma \ref{GT_ineq_gen_before}.

\hfill \qed

Note that the inequality (\ref{gen-GT-ineq}) can be regarded as a kind of one-parameter extension of the Golden-Thompson inequality for
positive semidefinite matrices $X$ and $Y$.

\section*{Ackowledgements}
We would like to thank the anonymous reviewer for providing valuable comments to improve the manuscript.
The first author (S.F.) was supported in part by the Japanese Ministry of Education, Science, Sports and Culture, Grant-in-Aid for
Encouragement of Young Scientists (B), 20740067

\end{document}